\documentclass[10pt,letterpaper,leqno]{article.altered}

\usepackage{amsmath}
\usepackage{amssymb}
\usepackage{amsthm}
\usepackage{accents}
\usepackage{bbding}
\usepackage{pifont}
\usepackage{booktabs}
\usepackage{ifthen}
\usepackage{enumerate.mod}
\usepackage{calc}
\usepackage[author-year]{amsrefs.altered}

\newtheorem{bigthm}{\large Theorem}
\newtheorem{thm}{Theorem}[section]
\newtheorem{lem}[thm]{Lemma}

\newtheorem{claim}[thm]{Claim}

\newtheorem{cor}[thm]{Corollary}

\newtheorem{question}{Question}

\theoremstyle{definition}
\newtheorem{defn}[thm]{Definition}

\newtheorem{notn}[thm]{Notation}

\theoremstyle{remark}
\newtheorem{remark}[thm]{Remark}

\newcommand{\ch}{\mathrm{CH}}

\newcommand{\ma}{\mathrm{MA}}
\newcommand{\zfc}{\mathrm{ZFC}}

\newcommand{\zf}{{\mathrm{ZF}}}

\newcommand{\pfa}{\mathrm{PFA}}
\newcommand{\rvm}{\mathrm{RVM}}

\newcommand{\sh}{\textrm{\rm SH}}

\newcommand{\pstar}{(*)}

\newcommand{\spr}{\bigvee}
\newcommand{\ifm}{\bigwedge}
\newcommand{\dom}{\operatorname{dom}}

\newcommand{\Fin}{\operatorname{Fin}}

\newcommand{\power}{\P}

\newcommand{\lbrak}{\bigl\|}
\newcommand{\rbrak}{\bigr\|}
\newcommand{\lmeas}{\mu\bigl(\lbrak}
\newcommand{\rmeas}{\rbrak\bigr)}

\newcommand{\inv}{^{-1}}

\newcommand{\otp}{\operatorname{otp}}

\newcommand{\lands}{\mathrel{\land}}
\newcommand{\lors}{\mathrel{\lor}}
\renewcommand{\div}{\mathbin{/}}
\newcommand{\cond}{\mathbin{|}}
\newcommand{\bigcond}{\mathbin{\big|}}

\newcommand{\Biggcond}{\mathbin{\Bigg|}}
\newcommand{\diff}{\bigtriangleup}

\newcommand{\pfn}{\dashrightarrow}
\newcommand{\elspc}[1]{\ell^{#1}}

\newcommand{\subseteqfnt}{\subseteq^*}

\newcommand{\supseteqfnt}{\mathrel{\supseteq^*}}

\newcommand{\Iff}{\espc\mathrm{iff}\espc}
\newcommand{\impls}{\qquad\text{implies}\qquad}
\newcommand{\AND}{\espc\mathrm{and}\espc}

\newcommand{\subseteqfntae}{\subseteqfnt_{\ae}}

\renewcommand{\b}{\mathfrak b}
\renewcommand{\d}{\mathfrak d}

\renewcommand{\t}{\mathfrak t}

\newcommand{\add}{\operatorname{add}}

\newcommand{\cantorcube}[1]{\{0,1\}^{\hsp#1}}
\newcommand{\two}{\{0,1\}}

\newcommand{\pN}{\power(\N)}
\newcommand{\reals}{\mathbb R}

\newcommand{\N}{\mathbb N}
\newcommand{\integers}{\mathbb Z}

\newcommand{\nulls}{\mathcal N}

\newcommand{\random}{\R}

\DeclareFontFamily{U}{cmsy}{}
\DeclareFontShape{U}{cmsy}{m}{n}{<12> sfixed * [10] cmsy10 
<10> <9> <8> <7> <6> <5> sfixed * [10] cmsy10}{}
\DeclareSymbolFont{customtwo}{U}{cmsy}{m}{n} 
\DeclareMathSymbol{\sctn}{\mathord}{customtwo}{"78}

\DeclareFontFamily{U}{cmmi}{}
\DeclareFontShape{U}{cmmi}{m}{n}{<20> sfixed * [11] cmmib10 <12> sfixed * [10]
cmmi10 <10> <9> <8> sfixed * [6] cmmi6 <5> <6> <7> sfixed * [5] cmmi5}{}
\DeclareSymbolFont{custom}{U}{cmmi}{m}{n}
\DeclareMathSymbol{\rharpoon}{\mathord}{custom}{"2A}
\newlength{\widt}
\newlength{\widttwo}
\newlength{\hgt}

\newcommand{\oone}{{\omega_1}}
\newcommand{\vomega}{\varOmega}

\newcommand{\<}{\langle}
\renewcommand{\>}{\rangle}

\newcommand{\espc}{\quad}

\renewcommand{\ae}{\mathrm{ae}}

\newcommand{\ulc}{\ulcorner}
\newcommand{\urc}{\urcorner}

\newcommand{\fsto}{\ding{75}}

\newcommand{\A}{\mathcal A}
\newcommand{\B}{\mathcal B}

\newcommand{\F}{\mathcal F}
\newcommand{\ideal}{\mathcal I}
\newcommand{\J}{\mathcal J}

\newcommand{\G}{\mathcal G}
\renewcommand{\H}{\mathcal H}

\renewcommand{\O}{\mathcal O}
\renewcommand{\P}{\Pcal}
\newcommand{\Pcal}{\mathcal P}
\newcommand{\Q}{\mathcal Q}
\newcommand{\R}{\mathcal R}
\newcommand{\Scal}{\mathcal S}

\newcommand{\hsp}{\mspace{1.5mu}}

\newcommand{\spc}{\,\,\,}

\newcounter{saveenumi}
\newcommand{\save}{\setcounter{saveenumi}{\value{enumi}}}
\newcommand{\restore}{\setcounter{enumi}{\value{saveenumi}}}

\newcommand{\tu}{\textup}

\DeclareFontFamily{T1}{pag}{}
\DeclareFontShape{T1}{pag}{m}{n}{<-> s * [0.95] phvr8t}{} 
\DeclareSymbolFont{special}{T1}{pag}{m}{n}
\DeclareMathSymbol{\nicecomplement}{\mathord}{special}{'143}

\renewcommand{\complement}{{\hspace{0.333pt}\nicecomplement}}

\newcommand{\texorpdfstring}[2]{#1}

\begin{document}

\title{\textbf{RANDOM GAPS}}
\author{\textsc{James Hirschorn}}
\date{}
\maketitle

\renewcommand{\thefootnote}{}

 \footnotetext{\textit{Date.} May 19, 2007.}
 \footnotetext{This research was primarily supported by Lise Meitner Fellowship, 
 Fonds zur F\"orderung der wissenschaftlichen Forschung, Project No.~M749-N05;
 the first version was completed on October 15, 2003,  
 with partial support of  Consorcio Centro de Investigaci\'on Matem\'atica, Spanish
 Government grant No.~SB2002-0099. Revisions were made
 with the support of Japanese Society for the Promotion of Science, 
 Project No.~P04301.}

\renewcommand{\thefootnote}{\arabic{footnote}}


\begin{abstract}
It is proved that there exists 
an $(\oone,\oone)$ Souslin gap in the Boolean algebra
$(L^0(\nu)\div\Fin,\subseteqfntae)$ for every nonseparable measure $\nu$. 
Thus a Souslin, also known as destructible, 
$(\oone,\oone)$ gap in $\pN\div\Fin$ can always be constructed from uncountably many random reals. We explain how to obtain the corresponding conclusion from
the hypothesis that Lebesgue measure can be extended to all subsets of the real line ($\rvm$).
\end{abstract}

\section{Introduction}
\label{sec:introduction}

A \emph{pregap} in a Boolean algebra $(\B,\le)$ is an orthogonal pair $(A,B)$ of
subsets of~$\B$, i.e.~%
\begin{enumerate}[(j)]
\item $a\lands b=0$ for all $a\in A$ and $b\in B$,
\save
\end{enumerate}
and it is a \emph{gap} if additionally there is no element $c$ of $\B$ such that
\begin{enumerate}[(j)]
\restore
\item $a< c$ for all $a\in A$, and $b< -c$ for all $c\in B$.
\save
\end{enumerate}
Such an element $c$ is said to \emph{interpolate} the pregap.
A \emph{linear pregap} is a pregap $(A,B)$ where both $A$ and $B$ are linearly
ordered by $\le$, and for a pair of linear order types $(\varphi,\psi)$, a
$(\varphi,\psi)$ \emph{pregap} in a Boolean algebra $(\B,\le)$ is a linear pregap
$(A,B)$ where $\otp(A,\le)=\varphi$ and $\otp(B,\le)=\psi$. 
Thus $(A,B)$ is a $(\varphi,\psi)$ \emph{gap} 
if it is a $(\varphi,\psi)$ pregap for which no element of
$\B$ can be used to extend $(A,B)$ to a $(\varphi+1,\psi)$ pregap or a
$(\varphi,\psi+1)$ pregap.

We let $\N$ denote the set $\{0,1,2,\dots\}$ of nonnegative integers, and $\pN$ the
power set of $\N$ quasi ordered by $a\subseteqfnt b$ if
$a\setminus b$ is finite.  $\pN\div\Fin$ denotes the equivalence classes of $\pN$ 
modulo the equivalence relation of finite set difference, with the induced partial ordering
$[a]\subseteqfnt [b]$ if $a\subseteqfnt b$.
One can find results on gaps in the Boolean algebra $(\pN\div\Fin,\subseteqfnt)$
dating back to the second half of the 19th century, including 
a basic result~(\citeauthor{DBR}) appearing in 1873, and
Hadamard's Theorem~\ycite{Had} that
there are no $(\omega,\omega)$ gaps in $\pN\div\Fin$.
Indeed, one of the major achievements in early Set Theory was
Hausdorff's construction~\ycite{Ha} of an $(\oone,\oone)$ gap in $\pN\div\Fin$
(he actually first constructed in~\ycite{Ha1} an $(\oone,\oone)$ gap in a different structure, 
that has a simple translation to a gap in $\pN\div\Fin$).
See~\ocite{Sch}\linebreak for the history of gaps.

\renewcommand{\footnoterule}{\vspace*{-3pt}\noindent\rule{.4\columnwidth}{0.4pt}\vspace*{2.6pt}}

\subsection{Acknowledgements}

We are grateful to Andr\'es Caicedo for explaining us how to improve our previous account of 
the axiomatization of random forcing by $\rvm$ (Theorem~\ref{u-2}).
The author wishes to thank the referee for doing an outstanding job. 

\subsection{Destructibility}
\label{sec:destructibility}

While being a pregap in $\pN\div\Fin$ is absolute between any `reasonable'
models---e.g.~transitive models of some large enough fragment of $\zf$---, the
property of being a gap is not. For example, if $(A,B)$ is an
$(\oone,\oone)$ gap in $\pN\div\Fin$ and $\Q$ is a poset which
collapses $\aleph_1$, then by Hadamard's theorem, 
which generalizes to any limit ordinal of countable cofinality, 
forcing with $\Q$ must introduce an element of $\pN\div\Fin$ which
interpolates $(A,B)$ and thus renders it a non-gap. Avoiding this particular
example, an $(\oone,\oone)$ pregap $(A,B)$ in $\pN\div\Fin$ is called
\emph{destructible} if there is an $\aleph_1$ preserving poset which
forces that $(A,B)$ is not a gap. 

The statement ``all $(\oone,\oone)$ gaps in $\pN\div\Fin$ are indestructible'' is in
fact a Ramsey theoretic statement which is closely analogous to Souslin's
Hypothesis. This becomes clear when one considers the characterization
of destructibility in Theorem~\ref{l-2} below.  We credit Theorems~\ref{u-3}
and~\ref{l-2} to Kunen~\ycite{Ku} and Woodin~\ycite{W}, respectively. 
The elegant Ramsey theoretic presentation is due to Todor\v cevi\'c~\ycite{T5}. 
See e.g.~\ocites{T5,Sch,TF} for proofs. 

When working with pregaps in
$\pN\div\Fin$ one often works with representatives (i.e.~subsets of $\N$) of the
equivalence classes. With every pair of families $a_i,b_i\subseteq\N$ indexed by
$i\in I$ we associate a partition $[I]^2=K_0\cup K_1$ via
\begin{equation}
  \label{eq:7}
  \{i,j\}\in K_0\Iff (a_i\cap b_j)\cup(a_j\cap b_i)=\emptyset.
\end{equation}
In the case where (the equivalence classes of) 
$(a_i,b_i:i\in I)$ is a pregap in $\pN\div\Fin$,  
we may assume---in order to avoid trivialities and thereby obtain more concise
results---that the representatives have been chosen so that
\begin{equation}
  \label{eq:8}
  a_i\cap b_i=\emptyset\espc\text{for all $i\in I$}.\tag*{\fsto}
\end{equation}
The following theorem characterizes $(\oone,\oone)$ gaps in $\pN\div\Fin$.

\begin{thm}[Kunen]
\label{u-3}
For every $(\oone,\oone)$ pregap in $\pN\div\Fin$ with representatives chosen 
satisfying condition~\/\ref{eq:8} the following are equivalent\textup{:
\begin{enumerate}[(a)]
\item \textit{$(a_\alpha,b_\alpha:\alpha<\oone)$ is a gap.}
\item\label{item:4} \textit{There is no uncountable $0$-homogeneous subset of\/ $\oone$.}
\item\label{item:5} \textit{The poset of all finite $1$-homogeneous subsets of\/ $\oone$ has the ccc.}
\end{enumerate}}
\end{thm}

\begin{remark}
Only the  case $I=\oone$ was considered because in general the theorem is
false, although for $(\kappa,\kappa)$ pregaps with $\kappa$ is regular and
uncountable, it does generalize by replacing ``uncountable'' with ``cardinality
$\kappa$'' and ``ccc'' with ``\text{$\kappa$-cc}''. 
\end{remark}

By switching the colors $0$ and $1$ one obtains a characterization of
destructibility. 

\begin{thm}[Woodin]
\label{l-2}
For every $(\oone,\oone)$ pregap in $\pN\div\Fin$ with representatives
satisfying~\/\ref{eq:8} the following are equivalent\textup{:
\begin{enumerate}[(a)]
\item \textit{$(a_\alpha,b_\alpha:\alpha<\oone)$ is destructible.}
\item\label{item:3} \textit{There is no uncountable $1$-homogeneous subset of\/ $\oone$.}
\item\label{item:2} \textit{The poset of all finite $0$-homogeneous subsets of\/ $\oone$ has the ccc.}
\end{enumerate}}
\end{thm}

\noindent Considering condition~\eqref{item:4} of Theorem~\ref{u-3}, 
condition~\eqref{item:2} of Theorem~\ref{l-2} says that
there is a poset with the ccc forcing that $(a_\alpha,b_\alpha:\alpha<\oone)$ is not
a gap. A Souslin tree itself has the ccc and
forces an $\oone$-branch through itself, which is analogous to forcing an
interpolation of a gap.  And condition~\eqref{item:3} of Theorem~\ref{l-2} can be
viewed in analogy with the property that Souslin trees have no uncountable
antichains. With this in mind, 
we provisionally refer to destructible and \emph{Souslin} $(\oone,\oone)$
pregaps interchangeably, until further definitions are made
in Section~\ref{sec:nonlinear-gaps}. 

\begin{remark}
\label{r-2}
  It is a theorem of Kunen~\ycite{Ku} (see also e.g.~Scheepers~\ycite{Sch})
  that for cardinals $\kappa$ and $\lambda$,
  every $(\kappa,\lambda)$ pregap in $\pN\div\Fin$ 
  with either $\kappa\ne\oone$ or $\lambda\ne\oone$ can be
  interpolated by a ccc poset. E.g.~when $\kappa=\lambda\ne\oone$, 
  the corresponding condition Theorem~\ref{l-2}\eqref{item:2} 
  (i.e.~replace ``$\oone$'' with ``$\kappa$'') holds for any such pregap.
\end{remark}

The analogy goes further when one considers the influence of four additional set
theoretic axioms: the existence of a diamond sequence ($\diamondsuit$); 
the Continuum
Hypothesis ($\ch$); the principle $\pstar$ for ideals of countable subsets of
$\oone$, a consequence of $\pfa$ which is consistent with $\ch$ and entails many
combinatorial consequences of $\pfa$ (see~\ocites{rst:A,MR1809418}); 
and Martin's Axiom for $\aleph_1$ many dense subsets ($\ma_{\aleph_1}$). 
Jensen proved that $\diamondsuit$ implies the existence of a
Souslin tree, while Todor\v cevi\'c proved in~\cite{D} (he took credit for this in
private communication) that $\diamondsuit$ implies the existence of
a destructible $(\oone,\oone)$ gap. It is a theorem of Jensen~\cite{DJ}) 
that $\sh$ is consistent with $\ch$, 
while Abraham--Todor\v cevi\'c~\ycite{rst:A} proved that $\pstar$ implies
both $\sh$ and ``all $(\oone,\oone)$ gaps are indestructible'' (and in particular
``all $(\oone,\oone)$ gaps are indestructible'' is consistent with $\ch$). It is a
theorem of Solovay--Tennenbaum~\ycite{rst:ST} that $\ma_{\aleph_1}$ implies $\sh$,
while it a theorem of Kunen~\ycite{Ku} that $\ma_{\aleph_1}$ implies 
that there are no Souslin $(\oone,\oone)$ gaps in
$\pN\div\Fin$. These results are summarized in the
\textbf{trivial} column of Tables~\ref{tab:1} and~\ref{tab:2}.

\begin{table}
\caption{\emph{Souslin's Hypothesis}}
\begin{center}
\label{tab:1}
\begin{tabular}{rlll} 
  \toprule  & \multicolumn{3}{c}{Measurable algebra}  \\
  \cmidrule(l){2-4}
                  \large{Axiom} & \textbf{trivial} & \textbf{separable} & \textbf{nonseparable} \\
  \midrule         $\diamondsuit$ & False                  & False                    & False \\
                            $\ch$ & Undecided       & Undecided                  & Undecided \\
                               $\zfc$ & Undecided       & Undecided            & Undecided \\
                          $\pstar$ & True                  & True            & True \\
         $\ma_{\aleph_1}$ & True         & True                       & True \\
  \bottomrule
\end{tabular}
\end{center}
\end{table}

\begin{table}
\caption{\emph{All gaps are indestructible}}
\begin{center}
\label{tab:2}
\begin{tabular}{rlll}
  \toprule       & \multicolumn{3}{c}{Measurable algebra}  \\                 
  \cmidrule(l){2-4}   
                  \large{Axiom} & \textbf{trivial} & \textbf{separable} & \textbf{nonseparable} \\
  \midrule 
                   $\diamondsuit$ & False                  & False                    & False \\
                   $\ch$             & Undecided       & False                      & False \\
                   $\zfc$            & Undecided       & Undecided            & False \\
                   $\pstar$        & True                  & Undecided            & False \\
                   $\ma_{\aleph_1}$ & True         & True                       & False \\
  \bottomrule
\end{tabular}
\end{center}
\end{table}

The analogy can be carried still further by considering the influence of Cohen forcing. 
It is a theorem of Shelah~\ycite{rst:S} that Cohen forcing (i.e.~adding one Cohen real)
always produces a Souslin tree, 
while it is a theorem of Todor\v cevi\'c (see~\ocite{TF}) that Cohen forcing
produces a Souslin gap. 

This naturally leads us at once to consider the influence of the other fundamental 
forcing notion, random forcing: 
a separable measurable algebra (adding one random real), 
or more generally forcing with an arbitrary measurable algebra (possibly
nonseparable, adding uncountably many random reals). 

Let us recall here that a Boolean algebra $\B$ is \emph{measurable} if
there is a function $\mu:\B\to[0,\infty)$ such that $(\B,\mu)$ is a measure algebra. 
The measure $\mu$ has a naturally associated metric $d_\mu$\label{dmu} on $\B$, 
where the distance between $a$ and $b$ is $d_\mu(a,b)=\mu(a\diff b)$.
The measure algebra is \emph{\tu(non\tu)separable}
if the metric topology on $\B$ is (non)separable. 
Note that given a measurable $\B$ separability is
independent of the choice of measure $\mu$; 
this follows from Maharam's Theorem which is discussed below.

The influence on Souslin's
Hypothesis of forcing with some measurable algebra over a model satisfying any of
the above axioms, is known and summarized in Table~\ref{tab:1}. It is a theorem of
Laver~\ycite{rst:L} that under $\ma_{\aleph_1}$, forcing with any measurable algebra
preserves $\sh$. 
And it is a theorem of the author~\ycite{H2} that assuming~$\pstar$, forcing with any
measurable algebra preserves $\sh$. Note that the other rows follow from these
results because Souslin trees are preserved by forcing notions satisfying property~$K$
(i.e.~Knaster's chain condition), and in particular by any measurable algebra.

It is a result of the author~\ycite{MR1949702} that under $\ma_{\aleph_1}$, all gaps 
are indestructible in any forcing extension by a separable measurable
algebra. In~\ycite{MR2034309} the author 
proves that under $\ch$, adding one random real produces a
Souslin gap (the $\ch$ row and \textbf{separable} column of Table~\ref{tab:2});
and this is the first place where the analogy breaks down, proving that:
\emph{Souslin's Hypothesis is consistent 
with the existence of a Souslin $(\oone,\oone)$ gap in $\pN\div\Fin$}. 
The following relatively old question of Woodin has received the attention of a number 
of authors. 

\begin{question}[Woodin]\label{q:1}
Does $\ma_{\aleph_1}$ imply that all $(\oone,\oone)$ gaps in $\pN\div\Fin$ are
indestructible in any forcing extension by a measurable algebra\textup?
\end{question}

The main result of this paper, the
``False'' in the $\zfc$ row of the \textbf{nonseparable} column of
Table~\ref{tab:2}, answers Question~\ref{q:1} negatively and it completes the
picture by allowing us to fill in the remainder of Table~\ref{tab:2}:

\begin{bigthm}\label{bigthm:1}
Let $\random$ be a nonseparable measurable algebra. 
Then some condition in $\random^+$
forces that there exists a Souslin $(\oone,\oone)$ gap in $\pN\div\Fin$. 
\end{bigthm}

\subsection{Terminology}
\label{sec:terminology}

We have been calling $(\B,\le)$ a Boolean algebra, 
to indicate that $(\B,\le)$ is a partial ordering with minimum and maximum
elements $0$ and $1$, resp., such that every two elements $x,y\in\B$ 
have both an infimum $x\lands y$ and a supremum $x\lors y$ and also has a complement
$-x=\spr\{z\in\B:x\lands z=0\}$.\footnote{We use the notation $x-y$ to abbreviate $x\lands(-y)$.}
This agrees with the usual definition of a Boolean algebra as a structure 
of the form $(\B,\land,\lor,-,0,1)$, in that $(\B,\le)$, 
with $x\le y$ defined by $x-y=0$, 
has the above properties 
iff this structure satisfies the axioms of a Boolean algebra. 

Recall that a \emph{measure space} is a triple $(X,\Sigma,\nu)$ where
$\Sigma$ is a $\sigma$-algebra of subsets of $X$ 
consisting of the
\emph{$\nu$-measurable sets}, $\nu:\Sigma\to[0,\infty]$ is a function 
with $\nu(\emptyset)=0$ and
\begin{equation}
  \label{eq:34}
  \nu\left(\bigcup_{n=0}^\infty E_n\right)=\sum_{n=0}^\infty\nu(E_n)
\end{equation}
whenever $(E_n:n\in\N)$ is a sequence of pairwise disjoint $\nu$-measurable sets.
A \emph{measure algebra} is a pair $(\B,\mu)$ where $\B$ is a
$\sigma$-complete Boolean algebra and $\mu:\B\to[0,\infty]$ is a
function with $\mu(0)=0$ and $\mu\bigl(\spr_{n=0}^\infty
a_n\big)=\sum_{n=0}^\infty \mu(a_n)$ 
whenever $a_m\land a_n=0$ for all
$m\ne n$. A \emph{probability space} is a measure space $(X,\Sigma,\nu)$ 
with $\nu(X)=1$; similarly, a \emph{probability algebra} is a measure algebra
$(\random,\mu)$ with $\mu(1)=1$. A \emph{$\sigma$-finite} measure
space is a measure space with $\nu$-measurable sets $(E_n:n\in\N)$
such that $X=\bigcup_{n=0}^\infty E_n$ and $\nu(E_n)<\infty$ for all~$n$; 
a $\sigma$-finite measure algebra is defined analogously.

For a measure space $(X,\nu)$ we write
$\nulls_\nu$ for the ideal of all $A\subseteq X$ with
$\nu(A)=0$. These sets are called \emph{null} (or \emph{negligible}).
The \emph{measure algebra of a measure space} 
$(X,\Sigma,\nu)$ is the quotient 
$\Sigma\div\nulls_\nu$ with the well-defined measure $\mu([E])=\nu(E)$.

We recall some further definitions we shall need. 
The \emph{additivity} of a measure algebra
$\add(\mu)$ is the the smallest cardinal $\kappa$ 
for which there is a pairwise incompatible
family $\A\subseteq\random$ of cardinality $\kappa$ such that the equality
\begin{equation}
  \label{eq:78}
  \mu\left(\spr\A\right)=\spr_{a\in\A}\mu(a)
\end{equation}
fails. Thus $\add(\mu)\ge\aleph_1$. 
A similar definition can be made for a measure space. 
Note that the additivity of the ideal $\nulls_\nu$ satisfies
$\add(\nulls_\nu)=\add(\nu)$. For some cardinal $\kappa$, we say that
a measure $\nu$ is \emph{$\kappa$-additive} to indicate that
$\add(\nu)\ge\kappa$. 
A measurable set $E\subseteq X$ is called an \emph{atom} if it is not
null and every $F\subseteq E$ with $\mu(F)<\mu(E)$ is null; similarly,
for a nonzero element of a measure algebra.

\subsection{Applications}
\label{sec:applications}

\subsubsection{Souslin trees versus gaps}
\label{sec:souslin-trees-versus}

An immediate consequence is a simple construction\footnote{The point is that it is
  relatively difficult to construct a model of $\pstar$ and $\ch$, and this
  consistency result was originally proved by adding one random real to such a model.} of a
model satisfying $\sh$ which also has a destructible gap:

\begin{cor}
Adding uncountably many random reals to any model of $\zfc$ satisfying
$\ma_{\aleph_1}$ gives a model of \textup{``}there are no Souslin trees\textup{''} and \textup{``}there is a
Souslin $(\oone,\oone)$ gap in $(\pN\div\Fin,\subseteqfnt)$\textup{''}. 
\end{cor}

\subsubsection{Consequences of $\rvm$}
\label{sec:consequences-rvm}

Another corollary is that the classical hypothesis that the Lebesgue measure can
be extended to all subsets of the real line $\reals$, implies 
the existence of a destructible $(\oone,\oone)$ gap in $\pN\div\Fin$. 
This is an immediate consequence of
Theorem~\ref{bigthm:1} and known absoluteness results for forcing
extensions by a
large enough measure algebra, from a real-valued measurable cardinal.  
We thankfully acknowledge Stevo Todor\v cevi\'c for suggesting 
Corollary~\ref{o-2} (in May 2002).

This absoluteness is in fact an axiomatization of random forcing. 
We state a theorem (Theorems~\ref{u-2})
to this effect, without proof, some variation of which is folklore. 
It says that if a  statement of
reasonable complexity is forced to hold in the extension by a large enough
measurable algebra, then this statement is a consequence of the axiom that the
Lebesgue measure extends to all subsets of the real line. 
This axiom is often called
$\rvm$. 

We will not give the proof here, or even a full explanation of the terminology. 
For a complete proof, and also a variation of Theorem~\ref{u-2}, 
with other aspects of the
mathematical background explained, we refer the
reader to the supplement to this paper~\cite{Hir-rg}; and for
further reading we also suggest the paper~\cite{MR2267147}. Recall
that an \emph{atomlessly measurable cardinal $\kappa$} is  
an uncountable cardinal carrying an atomless $\kappa$-additive
probability measure with domain $\power(\kappa)$. 
It is a classical theorem of Ulam~\ycite{U} that 
$\rvm$ is equivalent to the existence of an atomlessly measurable
cardinal $\kappa$, and that $\kappa$ is larger than the least weakly
inaccessible cardinal but~$\kappa\le 2^{\aleph_0}$.

\begin{notn}
\label{notn:1}
For a set $X$, we write $(\random_{(X)},\mu_{(X)})$ for the measure algebra
of the measure space $\cantorcube X$ with its Haar probability
measure.
\end{notn}

A word on Maharam's Theorem is also in order here. It states roughly
that every $\sigma$-finite measure algebra
(cf.~\S\ref{sec:terminology}), up to an isomorphism, has a simple
decomposition into measure algebras of the form
$(\random_{(\theta)},\mu_{(\theta)})$, where $\theta$ is some cardinal. A measurable
algebra (cf.~\S\ref{sec:destructibility}) is \emph{homogeneous} in
the forcing sense iff it is isomorphic to $\random_{(\theta)}$ for
some cardinal $\theta$. The
\emph{Maharam type} of such a measurable  algebra is the 
cardinal $\theta$. By Maharam's Theorem, every $\sigma$-finite
measurable algebra $\random$ 
has a dense set of $z\in\random^+$ such that
$\random_z=\{x\in\nobreak\random:x\le z\}$ is homogeneous.
Note that a homogeneous measurable algebra is nonseparable iff 
its Maharam type is at least $\aleph_1$.
For more information on Maharam's
Theorem, see also~e.g.~\ocites{F3,Hegt}.

\begin{thm}
\label{u-2}
Suppose $\kappa$ is an atomlessly measurable cardinal.
If $\varphi(x,y)$ is a $\varPi_1$ formula and $a\in H_\kappa$ is a
parameter such that every homogeneous
measurable algebra of large enough Maharam type forces the statement
\begin{equation}
  \label{eq:19}
  \exists x\in H_{\check\kappa^+}\spc \varphi(x,\check a),
\end{equation}
then $V\models\ulc\exists x\,\varphi(x,a)\urc$.
\end{thm}

\begin{remark}
\label{r-6}
There is an important theorem of
Gitik and Shelah~\ycite{GS1} (see
also~\ocite{MR1234282}) giving a lower bound on the Maharam type of the measure
algebras associated with atomlessly measurable cardinals: \emph{Let $X$ be a set. 
If $\nu:\power(X)\to[0,1]$ is an atomless probability measure then 
the measure algebra $\power(X)\div\nulls_\nu$ has Maharam type at least
$\min\bigl\{\add(\nu)^{+\omega},2^{\add(\nu)}\bigr\}$}. (Moreover they proved
in~\ocite{GS} that the Maharam type is $2^{2^{\aleph_0}}$ in the case
$X=2^{\aleph_0}$ and $\add(\nu)=2^{\aleph_0}$.)
The ``large enough Maharam type'' in Theorem~\ref{u-2} is the Maharam
type of $\power(\kappa)\div\nulls_\nu$, and thus is at least
$\min\bigl\{\kappa^{+\omega},2^\kappa\bigr\}$. 
\end{remark}

\begin{cor}[$\zfc+\rvm$]
\label{o-2}
There exists a Souslin $(\oone,\oone)$ gap in $\pN\div\Fin$.
\end{cor}
\begin{proof}
Using Woodin's Theorem (Theorem~\ref{l-2}), 
it is routine to obtain a $\varPi_1$ formula $\varphi(x,y,a)$ so that
$\varphi(x,y,\oone)$ holds iff $(x,y)$ forms a Souslin
$(\oone,\oone)$ gap in $\pN\div\Fin$. 
By the assumption that  the Lebesgue measure on the real line can be
extended to a measure whose domain is all of $\P(\reals)$, there
exists an atomlessly measurable cardinal $\kappa$. 
Since $(\oone,\oone)$ gaps are objects of $H_{\aleph_2}$,
by Theorem~\ref{bigthm:1},
every homogeneous measure algebra of Maharam type at least $\aleph_1$
forces $\ulc\exists x,y\in H_{\check\kappa^+}
\spc \varphi(x,y,\check\omega_1)\urc$.
Therefore, from Theorem~\ref{u-2} we conclude that there exists
$(x,y)$ satisfying $\varphi(x,y,\oone)$, completing the proof.
\end{proof}

\subsection{Souslin gaps}

\label{sec:nonlinear-gaps}

If $(A,B)$ is a gap, not necessarily linear, where the cardinalities of $A$ and $B$
are both at most $\aleph_1$, then forcing with a poset 
which collapses $\aleph_1$ will still interpolate the gap. 
Hence, it makes sense to extend the definition 
of \emph{destructible} to include this situation. 

Broadening our scope to include arbitrary pregaps in $\pN\div\Fin$,
Theorem~\ref{l-2} is no longer true. In condition~\eqref{item:2}, 
the poset of finite $0$-homogeneous subsets will force an interpolation
of an uncountable subgap, but if the gap is nonlinear it will not necessarily
interpolate the whole gap. We would like to isolate this condition~\eqref{item:2},
as the Souslin property.

\begin{defn}
\label{d-4}
A pregap $(A,B)$ in $\pN\div\Fin$ is called \emph{Souslin} if every uncountable
family $\F$ of finite subsets of $A\times B$ (i.e.~finite subpregaps of $(A,B)$) 
has an $\aleph_1$ preserving forcing extension
where there exists $d\in\pN\div\Fin$ such that
\begin{equation}
  \label{eq:71}
  \{F\in\F:d\text{ interpolates }F\}\text{ is uncountable.}
\end{equation}
Note we need only consider families $\F$ of size $\aleph_1$, and that by a
$\Delta$-system argument we can assume $\F$ is pairwise disjoint. 
\end{defn}

This Souslin property of gaps can be characterized in terms of representative
subsets of $\N$. For this it will be convenient to reconsider the
requirement~\ref{eq:8}. Suppose $(a_i,b_i:i\in I)$ is an indexing of two sequences of
members of $\pN$. 
Instead of partitioning the subsets of cardinality two $[I]^2$ as in~\eqref{eq:7},
we can define partitions of the pairs $I^2=L^k_0\cup L^k_1$ and generalize 
by adding a parameter $k\in\N$, as in
\begin{equation}
  \label{eq:57}
  (i,j)\in L^k_0\Iff (a_i\cap b_j\setminus k)\cup(a_j\cap b_i\setminus k)=\emptyset,
\end{equation}
and we write $L_0$ and $L_1$ for $L_0^0$ and $L_1^0$, respectively.
Then the condition~\ref{eq:8} is equivalent to ``$L_0(i,i)$ for all $i\in I$'', 
and assuming this is satisfied, a subset of $I$ is
$K_0$-homogeneous iff it is $L_0$-homogeneous.
Notice that the quantification in the next theorem is significantly
different than in Theorems~\ref{u-3} and~\ref{l-2}.

\begin{thm}
  \label{u-6}
Let $(A,B)$ be a pregap in $\pN\div\Fin$. Then the following are equivalent\textup{:
  \begin{enumerate}[(a)]
  \item\label{item:12} \textit{$(A,B)$ is Souslin.}
  \item\label{item:11} \textit{For every choice $(a_i,b_i:i\in I)$ of
      representatives the associated poset 
      of all finite $L_0$-homogeneous subsets of $I$ has the ccc.}
  \item\label{item:9} \textit{There exists an indexing $(a_i,b_i:i\in I)$ of
      representatives such that for every~$k$ the associated poset of all finite
      $L_0^k$-homogeneous subsets of $I$ has the ccc.} 
  \end{enumerate}}
\end{thm}
\begin{proof}
\eqref{item:12}$\,\to\,$\eqref{item:11}: Suppose that $(A,B)$ is Souslin and $(a_i,b_i:i\in I)$ is a choice of
representatives.
Let $J_\alpha\in\Fin_I$ ($\alpha<\oone$) be
a family of $0$-homogeneous finite subsets of $I$. 
By going to an uncountable subset
we can assume that they are all the same size, say
$J_\alpha=\{i^\alpha_0,\dots,i^\alpha_{n-1}\}$ for all $\alpha$. 
Now go to an $\aleph_1$
preserving forcing extension with $d\subseteq\N$ such that
\begin{equation}
  \label{eq:73}
  X=\bigl\{\alpha<\oone:d\text{ interpolates }
  \{(a_i,b_i):i\in J_\alpha\}\bigr\}\text{ is uncountable}.
\end{equation}
By going to an uncountable subset of $X$ we can assume there is a $k\in\N$ such that
$a_i\setminus k\subseteq d$ and $b_i\setminus k\subseteq d^\complement$ for all
$i\in J_\alpha$ for all $\alpha\in X$. Finally, we can pick $\alpha\ne\beta$
in $X$ such that both $\{(a_i,b_i):i\in J_\alpha\}$ and
$\{(a_i,b_i):i\in J_\beta\}$ have the same trace on $k$. Then
$(J_\alpha\cup J_\beta)^2\subseteq K_0$ because for all $l,m=0,\dots,n-1$,
\begin{equation}
  \label{eq:74}
  \begin{split}
  a_{i^\alpha_l}\cap b_{i^\beta_m}&=(a_{i^\alpha_l}\cap b_{i^\beta_m} \cap k)
  \cup(a_{i^\alpha_l}\cap b_{i^\beta_m}\setminus k)\\
  &=a_{i^\alpha_l}\cap b_{i^\beta_m}\cap k\cup\emptyset\\
  &= a_{i^\alpha_l}\cap b_{i^\alpha_m}\cap k\\
  &=\emptyset,
  \end{split}
\end{equation}
as $a_{i^\alpha_l}\cap b_{i^\alpha_m}=\emptyset$ for all $l,m$ by homogeneity.

\eqref{item:11}$\,\to\,$\eqref{item:9}: Take \emph{any} indexing $(a_i,b_i:i\in I)$
of representatives. Suppose $\F$ is an uncountable family of finite
$L_0^k$-homogeneous subsets of $I$ for some fixed $k\in\N$. Applying
condition~\eqref{item:11} to the indexing $(a_i\setminus k,b_i\setminus k:i\in I)$
of representatives, since each $F\in\F$ is $0$-homogeneous for the associated
partition, it follows that there exist $F\ne G$ in $\F$ such that 
$(F\cup G)^2\subseteq L_0^k$. 

\eqref{item:9}$\,\to\,$\eqref{item:12}: Suppose that $(a_i,b_i:i\in I)$ is an indexing of representatives satisfying the
hypothesis of~\eqref{item:9}. 
For each $F\in\F$ choose a finite $J_F\subseteq I$ such that 
$F\subseteq\{a_i:i\in J_F\}\times\{b_i:i\in J_F\}$, and then find $k_F\in\N$ large
enough so that
\begin{equation}
  \label{eq:82}
  J_F^2\subseteq L_0^{k_F}.
\end{equation}
Pick $k\in\N$ large enough so that $\F_0=\{F\in\F:k_F=k\}$ is uncountable.
It follows from condition~\eqref{item:9}
that the poset of all finite $\H\subseteq\F_0$ such that $\bigcup_{F\in\H}J_F$
is $L^k_0$-homogeneous has the ccc, and therefore it can have only countably many atoms
and thus forces an uncountable subset $\G\subseteq\F_0$
each member of which is interpolated by $\bigcup_{F\in\G}\bigcup_{i\in
  J_F}a_i$. 
\end{proof}

Pregaps have been called \emph{Luzin} in the literature (e.g.~\ocite{MR2048515}) when there
exists an uncountable $K_1$-homogeneous subset of $I$ for some indexing of 
representatives satisfying\/~\ref{eq:8}. By
Theorem~\ref{u-6}~\eqref{item:12}$\,\to\,$\eqref{item:11}, 
every Souslin pregap in $\pN\div\Fin$ is non-Luzin, 
and thus for pregaps of size at most~$\aleph_1$
\begin{equation}
  \label{eq:75}
  \text{destructible}\to\text{Souslin}\to\text{non-Luzin.}
\end{equation}
In the realm of linear gaps, Souslin gaps are precisely the non-Luzin ones, and for
linear gaps of size $\aleph_1$ all three notions coincide by
Theorem~\ref{l-2}. It may be worth revisiting this concept in a separate article.

\subsection{Gaps in \texorpdfstring{${(L^0(\nu)\div\Fin,\subseteqfntae)}$}{L(m)/Fin,)}}
\label{sec:random-gaps}

Let $(\random,\mu)$ be the measure algebra of some $\sigma$-finite 
measure space $(X,\Sigma,\nu)$. 
In considering $\random$-names for gaps in $\pN\div\Fin$,
we are led to consider gaps in another Boolean algebra. Recall that a \emph{random
  variable on $X$ with codomain~$S$}, where $S$ is some topological space, is an
almost everywhere defined function $f$, measurable for the completion of $\mu$, 
taking values in $S$; in other words, $f$ is
defined on a conegligible subset of $X$ and there
is a conegligible $C\subseteq X$ such that for every open $U\subseteq
S$, $f\inv[U]\cap C\in\Sigma$. 
The family of all equivalence classes of random
variables on $X$ with codomain $\reals$ over the equivalence relation $f(x)=g(x)$
almost everywhere, is the standard space $L^0(\nu)$ from functional analysis. 
We are interested here in taking $\pN$ to be our reals; hence, we define
$L^0(\nu)\div\Fin$ to be the family of all 
equivalence classes of random variables on $X$ with codomain $\pN$ 
modulo the relation
$f(x)\diff g(x)$ is finite almost everywhere, and order $L^0(\nu)\div\Fin$ by defining
$[f]\subseteqfntae[g]$ if $f(x)\subseteqfnt g(x)$ almost everywhere. This clearly
defines a Boolean algebra. By a \emph{random gap} we mean a (pre)gap in a Boolean
algebra of this form.

We have seen that Maharam's Theorem says that for purposes of forcing, 
we can restrict our attention to measure algebras of the form
$(\random_{(\theta)},\mu_{(\theta)})$ for some cardinal $\theta$. And it is well
known (e.g.~\ocite{MR0218233}) that there is a direct correspondence between
$\random_{(\theta)}$-names $\dot a$ for a subset of $\N$ 
and random variables $f$ on the underlying measure space
$\cantorcube\theta$ with codomain $\pN$, via
\begin{equation}
  \label{eq:32}
  \lbrak f(\dot r)=\dot a\rbrak=1,
\end{equation}
where $\dot r$ is a name for the generic object in $\cantorcube\theta$.
Of course $f$ must be interpreted correctly in the forcing extension to make sense
of~\eqref{eq:32}, but the point is that every measurable function 
$f:\cantorcube\theta\to\pN$ is equal almost everywhere to a Baire
function\footnote{We mean Baire measurable function, where the Baire sets are
  members of the smallest $\sigma$-algebra for which every continuous function is
  measurable. This is the $\sigma$-algebra generated by the clopen sets in the case
  of  a Cantor cube $\cantorcube X$.}
$g:\cantorcube\theta\to\pN$, and every such Baire
function can be coded by a countable sequence of ordinals
so that every suitable model containing this sequence has
a correct interpretation of $g$ as a Baire function from $\cantorcube\theta$ into $\pN$. 
Henceforth, we will dot random variables (e.g.~$\dot a$) to emphasize this
correspondence. 

For two of these random variables $\dot a$ and $\dot b$, clearly $\lbrak \dot
a(\dot r)\subseteqfnt \dot b(\dot r)\rbrak=1$ iff $[\dot a]\subseteqfntae[\dot b]$ in
$L^0(\mu_{(\theta)})\div\Fin$. 
Thus a pregap in $(L^0(\mu_{(\theta)})\div\Fin,\subseteqfntae)$ is the same thing as an
$\random_{(\theta)}$-name for a pregap in $(\pN\div\Fin,\subseteqfnt)$. 
And a pregap in $L^0(\mu_{(\theta)})\div\Fin$ is a gap iff when viewed as an
$\random_{(\theta)}$-name for a pregap it is a gap in $\pN\div\Fin$ 
with positive probability.
We also would like to extend the correspondence to the notion of Souslin gaps, 
keeping in mind Theorem~\ref{u-6}. 

\begin{defn}
A pregap $(A,B)$ in $L^0(\nu)\div\Fin$ is \emph{Souslin} if there is an enumeration
$(\dot a_\alpha,\dot b_\alpha:\alpha<\oone)$ of representatives, so that for every sequence
$(E_\xi,\varGamma_\xi:\xi<\nolinebreak\oone)$ and every $k\in\N$, where
$E_\xi\in\Sigma$ has positive finite measure, $\varGamma_\xi\subseteq\oone$ is
finite, and
\begin{equation}
  \label{eq:33}
  \Pr\Biggl[\bigcup_{\substack{\alpha,\beta\in\varGamma_\xi}}
  \dot a_\alpha\cap\dot b_\beta\setminus k=\emptyset\Biggcond E_\xi\Biggr]=1
\espc\text{for all $\xi$}
\end{equation}
(equivalently, $\bigcup_{\alpha,\beta\in\varGamma_\xi}
  \dot a_\alpha(x)\cap\dot b_\beta(x)\setminus k=\emptyset$ almost everywhere in $E_\xi$),
there exist $\xi\ne\eta$ such that $\nu(E_\xi\cap E_\eta)\ne0$ and 
\begin{equation}
  \label{eq:36}
  \Pr\Biggl[
  \bigcup_{\alpha\in\varGamma_\xi}
  \bigcup_{\beta \in\varGamma_\eta}
  (\dot a_\alpha\cap\dot b_\beta\setminus k)
  \cup(\dot a_\beta\cap\dot b_\alpha\setminus k)
  =\emptyset\Biggcond E_\xi\cap E_\eta\Biggr]\ne0.
\end{equation}
\end{defn}

Note that in the case where $\nu=\mu_{(\theta)}$, when identifying random variables
$\dot a$ on $\cantorcube\theta$ with codomain $\pN$ with the
$\random_{(\theta)}$-name $\dot a(\dot r)$ for a subset of $\N$, then replacing each
$E_\xi$ with $x_\xi\in\random^+$, condition~\eqref{eq:33} becomes
\begin{equation}
  \label{eq:52}
  x_\xi\lands\spr_{\alpha,\beta\in\varGamma_\xi }
  \lbrak \dot a_\alpha\cap\dot b_\beta\setminus k\ne\emptyset\rbrak=0\espc\text{for all $\xi$},
\end{equation}
and condition~\eqref{eq:36} becomes, with $x_\xi\lands x_\eta\ne0$,
\begin{equation}
  \label{eq:54}
  x_\xi\lands x_\eta-\spr_{\alpha\in\varGamma_\xi}\spr_{\beta\in\varGamma_\eta}
  \lbrak (\dot a_\alpha\cap \dot b_\beta\setminus k)
  \cup(\dot a_\beta\cap \dot b_\alpha\setminus k)
  \ne\emptyset\rbrak\ne0.
\end{equation}
Thus, by Theorem~\ref{u-6}, a pregap in $L^0(\mu_{(\theta)})$ is Souslin 
iff it is forced to be a Souslin
pregap with probability one when viewed 
as an $\random_{(\theta)}$-name.

For a measure space $(X,\Sigma,\nu)$, $\nu$ is called nonseparable if the
measure algebra of the measure space is nonseparable. Thus we obtain the following
equivalent formulation of Theorem~\ref{bigthm:1}.

\begin{bigthm}
\label{bigthm:2}
There exists an $(\oone,\oone)$ Souslin gap in $(L^0(\nu)\div\Fin,\subseteqfntae)$ 
whenever $\nu$ is a nonseparable $\sigma$-finite measure. 
\end{bigthm}

Let us note that gaps in $L^0(\nu)\div\Fin$ are a generalization of gaps in
$\pN\div\Fin$, because if $(X,\Sigma,\nu)$ is a trivial measure space with one
element of measure one, then $(L^0(\nu)\div\Fin,\subseteqfntae)$ and
$(\pN\div\Fin,\subseteqfnt)$ are isomorphic. Finally, we notice that the existence
of an $(\oone,\oone)$ gap in $L^0(\nu)\div\Fin$ for any $\sigma$-finite measure $\nu$, follows from
Hausdorff's Theorem on the existence of these gaps in $\pN\div\Fin$ applied in the
appropriate random forcing extension.

\subsubsection{Destructible gaps in \texorpdfstring{$L^0(\nu)\div\Fin$}{L(m)/Fin}}
\label{sec:destr-gaps}

The fact that there is no gap $(A,B)$ in $\pN\div\Fin$ with both $A$ and $B$
countable generalizes to gaps in $L^0(\nu)\div\Fin$ for any $\sigma$-finite measure
$\nu$. Therefore, an $(\oone,\oone)$ gap in $L^0(\nu)\div\Fin$ is interpolated by
collapsing $\aleph_1$. Thus we can define a destructible $(\oone,\oone)$ gap in
$L^0(\nu)\div\Fin$ as one which can be interpolated by a forcing which preserves
$\aleph_1$. Note that, in $L^0(\nu)\div\Fin$, destructible $(\oone,\oone)$ gaps and
Souslin $(\oone,\oone)$ gaps are two different things! 
Indeed in proving Theorem~\ref{bigthm:2} we shall construct a
Souslin gap which is indestructible. On the other hand, in~\ocite{MR2034309}, an
$(\oone,\oone)$ gap in $L^0(\nu)\div\Fin$ with $\nu$ separable, is constructed with
the aid of the Continuum Hypothesis, and this gap is both destructible and
Souslin; moreover, this gap can be destroyed by a poset with the property $K$, which
cannot happen to a gap in $\pN\div\Fin$. It is not hard to prove every destructible
$(\oone,\oone)$ gap in $L^0(\nu)\div\Fin$ must be Souslin.

\section{Measure theoretic characterizations}
\label{sec:meas-theor-char}

There are some simple measure theoretic characterizations, 
or more precisely necessary and/or sufficient conditions, of when 
an $\random$-name for a pregap names a Souslin gap and related phenomenon. 
This is our explanation for the complex behavior observed in the interactions 
between gaps in $\pN\div\Fin$ and random forcing.
For example, in~\ocite{MR2034309} it is shown that for a 
separable probability algebra $(\random,\nu)$,
an $\random$-name $(\dot a_\alpha,\dot b_\alpha:\alpha<\nolinebreak\oone)$ 
for an
pregap in $\pN\div\Fin$ is forced to be Souslin with positive probability, if
\begin{equation}
  \label{eq:37}
  \lmeas\dot a_\alpha\cap\dot b_\beta\ne\emptyset\rmeas<\frac12
  \espc\text{for all $\alpha,\beta<\oone$},
\end{equation}
(or less than some other constant below 1).
The statement~\eqref{eq:37} is far
simpler than the definition of destructibility or its Ramsey theoretic
characterization in Theorem~\ref{l-2}. 
Note that unlike Theorem~\ref{l-2}, the characterization 
in equation~\eqref{eq:37} is by no means equivalent to destructibility; 
in the case of the trivial
measure algebra it implies that the pregap is a non-gap.

We have also obtained in~\ocite{H7} a simple necessary condition 
for an $\random$-name for a pregap
to satisfy the \emph{Hausdorff property}, the strengthening of indestructibility to
$\{\alpha<\beta:a_\alpha\cap b_\beta\subseteq k\}$ is finite for all $k\in\N$, for
all $\beta<\oone$, which was satisfied by Hausdorff's original construction of an
$(\oone,\oone)$ gap. It states that if the equation~\eqref{eq:10} holds for
some $h\in c_0$ (i.e.~$\lim_{n\to\infty}h(n)=0$), 
then we do not have a name for a gap satisfying Hausdorff's
condition. In fact, in~\ocite{H7} an $\random$-name (with $\random$ nonseparable) 
for an indestructible gap not satisfying Hausdorff's condition is constructed by
satisfying~\eqref{eq:10} for an $h\in\elspc1$ (i.e.~$\sum_{n=0}^\infty h(n)<\infty$).
Therefore, the characterization of the Souslin property in Theorem~\ref{l-1} 
is not valid with the hypothesis~\eqref{eq:10} alone.

We were unable to find a sufficient condition for the Souslin property 
with an arbitrary measurable algebra, 
purely by putting constraints on the measures of various events.
This lead us to probabilistic considerations in our goal to obtain a
characterization of Souslin gaps in  $L^0(\nu)\div\Fin$ (this is
achieved in Theorem~\ref{l-1} below). 

Recall that a family $\{\Scal_i:i\in I\}$ of subsets of some measure
algebra $(\random,\mu)$ is \emph{stochastically independent} 
if $\mu(x\lands y)=\mu(x)\mu(y)$ for all
$x\in \Scal_i$ and $y\in\Scal_j$, for all $i\ne j$ in $I$. And two
elements $x,y\in\random$ are called stochastically independent when
$\bigl\{\{x\},\{y\}\bigr\}$ is a stochastically independent family. 

The following basic result of probability theory is not used for the
characterization of a Souslin gap, but is needed later on in the
actual construction of an $(\oone,\oone)$ Souslin gap in
$L^0(\nu)\div\Fin$ (a proof can be found in~\cite{Hir-rg}).

\begin{lem}
\label{l-5}
Let $\bigl\{\{x_i,y_i\}:i\in I\bigr\}$ be a stochastically independent family such
that $x_i\lands y_i=0$ for all $i\in I$. Then
for all $A,B\subseteq I$, $\mu\bigl(\spr_{i\in A}x_i\lands\spr_{i\in
  B}y_i\bigr)\le\mu\bigl(\spr_{i\in A}x_i\bigr)\mu\bigl(\spr_{i\in B}y_i\bigr)$
with equality iff either $A\cap B=\emptyset$ or $x_i=1$ or $y_i=1$ for some $i\in A\cup B$.
In other words,
$\spr_{i\in A}x_i$ is \emph{unfavourable} for $\spr_{i\in B}y_i$.
\end{lem}

Recall that for $A\subseteq X$, 
a subset $S$ of the Cantor cube $\cantorcube X$ is called
\emph{determined by coordinates in $A$}\label{p:1} 
or \emph{$A$-determined} 
if for all $y,z\in\cantorcube X$,
\begin{equation}
  \label{eq:44}
  y\restriction A=z\restriction A\impls y\in S\Iff z\in S.
\end{equation}
And $x\in\random_{(X)}$ is  \emph{determined by coordinates in $A$} or
\emph{$A$-determined} if $x$ is a member of the natural identification of
$\random_{(A)}$ with a subalgebra of $\random_{(X)}$ (cf.~Notation~\ref{notn:1}). 
Thus $x\in\random_{(X)}$ is $A$-determined 
iff it has an $A$-determined representative $S\subseteq\cantorcube X$. 
The set of \emph{coordinates determining $x$} is the minimum $A\subseteq X$
which determines~$x$. Such an $A$ always exists; 
indeed it is the set of all $\xi\in X$ such that
\begin{equation}
  \label{eq:161}
  \pi_{\<\xi,0\>}(x\lands \<\xi,0\>)\ne \pi_{\<\xi,1\>}(x\lands\<\xi,1\>),
\end{equation}
where $\pi_{\<\xi,i\>}:(\random_{(X)})_{\<\xi,i\>}
\to\random_{(X\setminus\{\xi\})}$ 
is the natural isomorphism as in\linebreak Corollary~\ref{o-1}, and $\<\xi,i\>$ 
denotes\label{pageref:<xi,0>} the
element of the measure algebra represented by $\{z\in\cantorcube X:z(\xi)=i\}$.
Recall that every member of $\random_{(X)}$ is represented by an $F_\sigma$ 
Baire set, and as such is determined by countably many coordinates.
The set of \emph{coordinates determining $\Scal$} 
where $\Scal\subseteq\random_{(X)}$
is the union of the coordinates that determine each member of $\Scal$.

We say that a collection $\{\Scal_i:i\in I\}$ of subsets of $\random_{(X)}$ 
is \emph{independently determined} if 
there is a family $\{A_i:i\in I\}$ of pairwise disjoint subsets of
$X$ such that $x$ is $A_i$-determined for all $x\in\Scal_i$,
i.e.~$\Scal_i\subseteq\random_{(A_i)}$, for all $i\in I$.
We also say that two elements $x,y\in\random_{(X)}$ are
\emph{independently determined} if the family
$\bigl\{\{x\},\{y\}\bigr\}$ is independently determined, and similarly
we can say that $x$ and $\Scal\subseteq\random_{(X)}$ are
\emph{independently determined} to indicate that
$\bigl\{\{x\},\Scal\bigr\}$ is independently determined.

Note that if $\F$ is an independently determined family of subsets then
$\F$ is stochastically independent; however, conversely, two stochastically
independent members of the measure algebra may fail to be disjointly determined.

For a finite partial function
$s:X\pfn\two$, we let $[s]=[s]_{(X)}\in\random_{(X)}$ 
be the equivalence class of the basic clopen set
\begin{equation}
\label{eq:31}
  \bigl\{z\in\cantorcube X:z\supseteq s\bigr\}.
\end{equation}
In some contexts, $[s]$ will denote the clopen set~\eqref{eq:31} instead. Let
$\Fin(X,\two)$ denote the collection of all finite partial functions from $X$ into
$\two$. 

\begin{lem}
\label{l-6}
Let $(\random,\mu)=(\random_{(\theta)},\mu_{(\theta)})$ 
for some infinite cardinal $\theta$,
and let\/ $t\in\Fin(\theta,\allowbreak \two)$. 
Then every $x\in\random_{[t]}$ has a unique
$\theta\setminus\dom(t)$-determined $y\in\random$ such that $x=y\lands[t]$.
\end{lem}


\begin{cor}
\label{o-1}
There is a unique measure algebra isomorphism
$\pi:\random_{[t]}\to\random_{(\theta\setminus(\dom(t))}$ such that
$\pi([s])=[s\setminus t]$ for all $s\in\Fin(\theta,\two)$ compatible with $t$. 
\end{cor}
\begin{proof}
Let $\pi(x)$ be the identification of the element of $\random$ given by
Lemma~\ref{l-6} with an element of $\random_{(\theta\setminus(\dom(t))}$.
\end{proof}

\noindent Maharam's Theorem gives an isomorphism between
$\random_{[t]}$ and $\random_{(\theta\setminus\dom(t))}$, 
but the point of the preceding corollary was to identify the
natural one.

\subsection{A chain condition for conditional probabilities}
\label{sec:chain-cond-cond}

There is a classical chain condition for probability algebras due to
Gillis~\ycite{G}, which
entails that every uncountable subset $\A$ of a probability algebra, 
in which every element has measure greater than some $\delta$, 
contains an uncountable subset $\B\subseteq\A$ where
$\mu(x\lands y)>\delta^2$ for every two elements $x$ and $y$ in $\B$. 
What we would
like here is an analogue for conditional probabilities:
\begin{equation}
  \label{eq:51}
  \mu(x\cond y)=\frac{\mu(x\lands y)}{\mu(y)},
\end{equation}
i.e.~for uncountable sequences $x_\alpha,y_\alpha\in\random^+$ $(\alpha<\oone)$ such
that $\mu(x_\alpha\cond y_\alpha)>\delta$ for all $\alpha$, there exists an
uncountable $X\subseteq\oone$ such that 
$\mu(x_\alpha\lands x_\beta\cond y_\alpha\lands y_\beta)>\delta^2$ for all $\alpha\in X$.
Although this is false as stated (cf.~\cite{Hir-rg}), 
if the $y_\beta$'s are of a simple enough form
then it becomes valid (Theorem~\ref{u-1}). 

\begin{thm}
\label{u-1}
Let $(\random,\mu)=(\random_{(\theta)},\mu_{(\theta)})$ for some cardinal $\theta$.
If $\A$ is an uncountable family of
pairs $(x,s)$ where $x\in\random$, $s\in\Fin(\theta,\two)$, 
and $\mu(x\cond[s])>\delta$ for all $(x,s)\in\A$, then there is an uncountable
$\B\subseteq\A$ such that $\mu(x\lands y\cond\allowbreak [s]\lands[t])>\delta^2$ 
for all $(x,s),(y,t)\in\B$.
\end{thm}
\begin{proof}
Let $(x_\alpha,s_\alpha)$ ($\alpha<\oone$) enumerate a subset of $\A$. 
Find an uncountable $X\subseteq\oone$ such that $\{\dom(s_\alpha):\alpha\in X\}$ 
forms a $\Delta$-system with root $\varGamma\subseteq\theta$, and
\begin{equation}
  \label{eq:162}
  \begin{aligned}
  |s_\alpha\setminus\varGamma|&=m\\
  s_\alpha\restriction\varGamma&=t
  \end{aligned}\espc\text{for all $\alpha\in X$}.
\end{equation}

Denote $\mu_{(\theta\setminus\varGamma)}$ by $\rho$, and write $\pi=\pi_{[t]}$ for
the isomorphism from Corollary~\ref{o-1}. Since the scalar correction of
$\mu$ to a probability measure on $\random_{[t]}$ is given  by $\mu(\cdot)\div\mu([t])$, 
we have $\rho\bigl(\pi(x)\bigr)=\frac{\mu(x)}{\mu([t])}$ for all $x\in\random_{[t]}$. 
However, by cancellation the scalar multiple of a measure has
the same conditional probabilities as the original measure. Hence, it suffices to
find an uncountable $X'\subseteq X$ such that
\begin{equation}
  \label{eq:90}
  \rho\bigl(\pi(x_\alpha)\lands \pi(x_\beta)\bigcond
  \pi([s_\alpha])\lands\pi([s_\beta])\bigr)>\delta^2
    \espc\text{for all $\alpha,\beta\in X'$.}
\end{equation}

Furthermore, the hypothesis translates to
$\rho(\pi(x_\alpha \lands [s_\alpha]))>\delta\rho(\pi([s_\alpha]))$ for all
$\alpha$. All of the $\rho(\pi([s_\alpha]))$'s for $\alpha\in X$ are equal to some
fixed $\sigma>0$---namely, $\sigma=2^{-m}$ by~\eqref{eq:162}.
Thus by Gillis' Theorem, there is an uncountable $X'\subseteq X$ such that
\begin{equation}
  \label{eq:17}
  \rho\bigl(\pi(x_\alpha\lands[s_\alpha])\lands\pi(x_\beta\lands[s_\beta])\bigr)
  >\delta^2\sigma^2\espc\text{for all $\alpha,\beta\in X'$.}
\end{equation}
But by the $\Delta$-system construct, $\pi([s_\alpha])=[s_\alpha\setminus t]$ is
determined independently of $\pi([s_\beta])=[s_\beta\setminus t]$ for all
$\alpha\ne\beta$ in $X$, and thus by stochastic independence,
$\rho(\pi([s_\alpha])\lands\pi([s_\beta]))=\sigma^2$. This establishes~\eqref{eq:90}. 
\end{proof}

\subsection{Continuous representatives}
\label{sec:cont-repr}

In~\ocite{MR1784706} we observed that every member 
of $L^0(\mu_{(\theta)})\div\Fin$ has a
continuous representative (actually the entire paper deals only with
$\theta=\omega$, but the proof of this fact applies to arbitrary $\theta$). Here we
make additional specifications on the representatives. The existence of continuous
representatives is essentially the fact that  the collection of
equivalence classes of clopen sets 
is dense in the metric topology on $\random_{(\theta)}$
(in other words, for every measurable set $A\subseteq\cantorcube\theta$ and every
$\varepsilon>0$ there exists a clopen $C\subseteq\cantorcube\theta$ such that
$\mu(A\diff C)<\varepsilon$).

\begin{lem}
\label{l-3}
For every pair of random variables $\dot a$,$\dot b$ on $\cantorcube\theta$ with
codomain $\pN$, and every $h:\N\to[0,1]$ with
\begin{equation}
  \label{eq:70}
  \lmeas n\in\dot a\cap\dot b\rmeas+h(n)>0\espc\textup{for all $n$,}
\end{equation}
there exist continuous functions $\dot c,\dot d:\cantorcube\theta\to\pN$ such that
\textup{\begin{enumerate}[(a)]
\item \textit{$\lbrak \dot c\cap\dot d=\emptyset\rbrak=1$,}
\save
\end{enumerate}}
\noindent and for all\/ $n$,
\textup{
\begin{enumerate}[(a)]
\restore
\item\label{item:1} 
  \textit{$\lbrak n\in\dot c\rbrak$ and $\lbrak n\in\dot d\rbrak$ are both
  determined by the set of coordinates which determines 
  $\bigl\{\lbrak n\in\dot a\rbrak,\lbrak n\in\dot b\rbrak\bigr\}$,}
\item \textit{$\mu\bigl(\lbrak n\in\dot a\diff\dot c\rbrak\bigr)
  <\lmeas n\in\dot a\cap\dot b\rmeas+h(n)$,}
\item \textit{$\mu\bigl(\lbrak n\in\dot b\diff\dot d\rbrak\bigr)
  <\lmeas n\in\dot a\cap\dot b\rmeas+h(n)$.}
\end{enumerate}}
\end{lem}
\begin{proof}
For each $n$: Let $A_n\subseteq\theta$ be the set of coordinates
determining $\bigl\{\lbrak n\in\nobreak\dot a\rbrak$, $\lbrak n\in\dot b\rbrak\bigr\}$. 
In case $A_n$ is finite, 
$\lbrak n\in\dot a\rbrak$ and $\lbrak n\in\dot b\rbrak$ both have
clopen representatives $C_n,D_n\subseteq\cantorcube\theta$,
respectively. Otherwise,
when $A_n$ is infinite, working in the space $\cantorcube{A_n}$, 
we can find a measurable $E_n\subseteq\lbrak n\in\dot a\cap\dot b\rbrak$ 
which is $A_n$-determined with measure
\begin{equation}
  \label{eq:46}
  \mu(E_n)=\frac{\lmeas n\in\dot a\cap\dot b\rmeas}2.
\end{equation}
Requirement~\eqref{eq:70} allows us to find $A_n$-determined 
clopen sets $C_n$ and $D'_n$ such that
\begin{align}
  \label{eq:47}
  \mu\Bigl(C_n \diff \Bigl(\lbrak n\in\dot a\rbrak
  \setminus\bigl(\lbrak n\in\dot b\rbrak\setminus E_n\bigr)\Bigr)\Bigr)
  &<\frac{\lmeas n\in\dot a\cap\dot b\rmeas+h(n)}6,\\
  \label{eq:48}
  \mu\Bigl(D'_n \diff \bigl(\lbrak n\in\dot b\rbrak\setminus E_n\bigr)\Bigr)
  &<\frac{\lmeas n\in\dot a\cap\dot b\rmeas+h(n)}6, 
\end{align}
i.e.~\eqref{eq:47} and~\eqref{eq:48} hold for any representatives of $\lbrak
n\in\dot a\rbrak$ and $\lbrak n\in\dot b\rbrak$. 
Then since $\lbrak n\in\dot a\rbrak \cap\bigl(\lbrak n\in\dot b\rbrak\setminus E_n\bigr)
=\lbrak n\in\dot a\cap\dot b\rbrak\setminus E_n$, which has measure  
$\lmeas n\in\dot a\cap\dot b\rmeas\div2$ by~\eqref{eq:46}, 
recalling the metric $d_\mu$ (cf.~page~\pageref{dmu}), 
the triangle inequality gives 
\begin{equation}
  \label{eq:49}
  \mu\bigl(C_n\diff\lbrak n\in\dot a\rbrak\bigr)
  <\frac{2\lmeas n\in\dot a\cap\dot b\rmeas+h(n)}3.
\end{equation}
And from~\eqref{eq:47} and~\eqref{eq:48} 
we have 
$\mu(C_n\cap D'_n)<\bigl[\lmeas n\in\dot a\cap\dot b\rmeas+h(n)\bigr]\div3$, 
and thus letting $D_n=D'_n\setminus C_n$ yields
\begin{equation}
  \label{eq:50}
  \begin{split}
    \mu\bigl(D_n\diff\lbrak n\in\dot b\rbrak\bigr)
    &\le\mu\bigl(D_n\diff \bigl(\lbrak n\in\dot b\rbrak\setminus E_n\bigr)\bigr)\\
    &\qquad+\mu\bigl(\bigl(\lbrak n\in\dot b\rbrak\setminus E_n\bigr)
    \diff \lbrak n\in\dot b\rbrak\bigr)\\
    &\le\mu(D_n\diff D'_n)
    +\mu\bigl(D'_n\diff\bigl(\lbrak n\in\dot b\rbrak\setminus E_n\bigr)\bigr)+\mu(E_n)\\
    &<\left(\frac13+\frac16+\frac12\right)\bigl[\lmeas n\in\dot a\cap\dot b\rmeas+h(n)\bigr]\\
      &=\lmeas n\in\dot a\cap \dot b\rmeas+h(n),
  \end{split}
\end{equation}
using the triangle inequality with $d_\mu$.

Now the functions $\dot c,\dot d:\cantorcube\theta\to\pN$ given by
$\dot c(x)=\{n:x\in C_n\}$ and $\dot d(x)=\{n:x\in D_n\}$ 
are continuous and as needed. 
\end{proof}

\begin{remark}
\label{r-1}
Note that if in Lemma~\ref{l-3}, 
$\sum_{n=0}^\infty\lmeas n\in\dot a\cap\dot b\rmeas<\infty$, 
then $\lbrak[\dot a]=[\dot c]\rbrak=1$, 
i.e.~$[\dot a]=[\dot c]$ in $L^0(\mu_{(\theta)})\div\Fin$, because
\begin{equation}
  \label{eq:9}
  \lbrak [\dot a]\ne[\dot c]\rbrak=\ifm_{k=0}^\infty\spr_{n=k}^\infty
  \lbrak n\in\dot a \diff \dot c\rbrak,
\end{equation}
and obviously $\lbrak[\dot b]=[\dot d]\rbrak=1$ too.
\end{remark}
\subsection{Characterization}
\label{sec:characterization}

The following Theorem~\ref{l-1} identifies some Souslin random gaps. 

\begin{bigthm}\label{l-1}
Let $\mu=\mu_{(\theta)}$ for some cardinal $\theta$.
Suppose $\{\dot a_\alpha:\alpha<\oone\}$ and $\{\dot b_\alpha:\alpha<\oone\}$ are two
families of representatives of members of $L^0(\mu)\div\Fin$. If
\begin{equation}
  \label{eq:56}
  \left\{\bigcup_{\alpha<\oone}\bigl\{\lbrak n\in\dot a_\alpha\rbrak
  \textup{, }\lbrak n\in\dot b_\alpha\rbrak\bigr\}:n\in\N\right\}
  \textup{ is independently determined,}
\end{equation}
and for some $h\in\elspc1$, 
\begin{equation}
  \label{eq:10}
  \lmeas n\in\dot a_\alpha\cap\dot b_\beta\rmeas\le h(n)
  \espc\textup{for all $\alpha,\beta<\oone$ and all $n\in\N$,}
\end{equation}
then $(\{[\dot a_\alpha]:\alpha<\oone\},\{[\dot b_\alpha]:\alpha<\oone\})$
forms a Souslin pregap.
\end{bigthm}

Note that condition~\eqref{eq:10} implies that $\{[\dot a_\alpha]:\alpha<\oone\}$
is orthogonal to $\{[\dot b_\alpha]:\alpha<\oone\}$ because $[\dot
a_\alpha]\lands[\dot b_\beta]=0$ in $L^0(\mu)\div\Fin$ is equivalent to 
$\lbrak\dot a_\alpha\cap\dot b_\beta$ is infinite$\rbrak=0$, and
\begin{equation}
  \label{eq:67}
  \lbrak\dot a_\alpha\cap\dot b_\beta\text{ is infinite}\rbrak
  =\ifm_{k=0}^\infty\spr_{n=k}^\infty\lbrak n\in\dot a_\alpha\cap\dot b_\beta\rbrak.
\end{equation}

Condition~\eqref{eq:56} is a very demanding requirement. If we consider an
$(\oone,\oone)$ gap in $L^0(\mu)\div\Fin$, i.e.~a linear gap, then we can relax this
to  $\bigl(\bigl\{\lbrak n\in\dot a_\alpha\rbrak$, $\lbrak n\in\dot
b_\alpha\rbrak\bigr\}:n\in\N\bigr)$  is independently determined sequence 
for every $\alpha$ separately.
We feel that such requirements are not very natural and that with the right
probability theory (concerning countable collections of random variables) a
natural characterization, say $\bigl\{\bigl\{\lbrak n\in\dot a_\alpha\rbrak$, $\lbrak
n\in\dot b_\alpha\rbrak\bigr\}:n\in\N\bigr\}$ is a stochastically independent family,
or even some further weakening, is obtainable. 
\subsubsection{Proof of Theorem~\ref{l-1}}
\label{sec:proof-theorem}

The proof begins by choosing a sequence
$(A^n:n\in\N)$ of pairwise disjoint subsets of $\theta$ such that 
\begin{equation}
  \label{eq:68}
  \lbrak n\in\dot a_\alpha\rbrak\text{ and }\lbrak n\in\dot b_\alpha\rbrak
  \text{ are both $A^n$-determined}\espc\text{for all $n$, for all $\alpha$}.
\end{equation}
Without loss of generality assume that $h>0$. 
Then for each $\alpha$, we apply Lemma~\ref{l-3} to the pair $(\dot a_\alpha,\dot b_\alpha)$, 
with $h_\alpha(n)=h(n)-\lmeas n\in\dot a_\alpha\cap\dot b_\alpha\rmeas$, 
to obtain a pair $(\dot c_\alpha,\dot d_\alpha)$ of continuous functions
as in the conclusion of the Lemma. Note that we have 
\begin{equation}
  \label{eq:12}
  \begin{split}
  \lmeas n\in\dot c_\alpha\cap \dot d_\beta\rmeas
  &\le\lmeas n\in\dot a_\alpha\cap \dot b_\beta\rmeas
  +\lmeas n\in\dot c_\alpha\setminus \dot a_\alpha\rmeas\\
  &\qquad\quad+\lmeas n\in\dot d_\beta\setminus \dot b_\beta\rmeas\\
  &< h(n)
  +\bigl(\lmeas n\in\dot a_\alpha\cap\dot b_\alpha\rmeas+h_\alpha(n)\bigr)\\
  &\qquad\quad
  +\bigl(\lmeas n\in\dot a_\beta\cap\dot b_\beta\rmeas+h_\beta(n)\bigr)\\
  &=3h(n)
  \end{split}
\end{equation}
for all $n$, for all $\alpha,\beta<\oone$.
For each $\alpha$, for each $n\in\N$, letting 
\begin{equation}
\label{eq:11}
y^n_\alpha=\lbrak n\in\dot c_\alpha\rbrak\AND z^n_\alpha=\lbrak n\in\dot d_\alpha\rbrak,
\end{equation}
by continuity we can write $y^n_\alpha=\spr_{s\in C^n_\alpha}[s]$ and
$z^n_\alpha=\spr_{s\in D^n_\alpha}[s]$ where
$C^n_\alpha,D^n_\alpha\subseteq\Fin(A^n,\two)$ are both finite, using
clause~\eqref{item:1} of Lemma~\ref{l-3}. Thus
\begin{equation}
  \label{eq:69}
  B_\alpha^n=\bigcup_{s\in C_\alpha^n \cup D_\alpha^n}\dom(s)\subseteq A^n
\end{equation}
is finite. 

By Remark~\ref{r-1}, 
it suffices to prove that $(\dot c_\alpha,\dot d_\alpha:\alpha<\oone)$ is Souslin. 
Suppose that $k\in\N$, $(x_\xi:\xi<\oone)$ is a sequence in $\random^+$ and
$(\varGamma_\xi:\xi<\oone)$ is a sequence of finite subsets of $\oone$ such
that
\begin{equation}
  \label{eq:39}
  x_\xi\lands\spr_{\alpha,\beta\in\varGamma_\xi}
  \lbrak\dot c_\alpha\cap\dot d_\beta\setminus k\ne\emptyset\rbrak=0\espc\text{for all $\xi$}.
\end{equation}
By going to an uncountable subsequence, 
we can assume that $|\varGamma_\xi|=m$ for all~$\xi$.

Let $\varepsilon>0$ be given. 
For each $\xi<\oone$, choose $s_\xi\in\Fin(\theta,\two)$ such that  
\begin{equation}
  \label{eq:59}
  \mu\bigl([s_\xi]\lands x_\xi\bigr)>\sqrt{1-\frac\varepsilon2}\cdot\mu([s_\xi]).
\end{equation}
Since the $A_n$'s are pairwise disjoint, if we choose a large enough 
$p_\xi\in\N$, then since $h\in\ell^1$ and by equation~\eqref{eq:69},
\begin{gather}
  \label{eq:21}
  \sum_{n=p_\xi}^\infty h(n)<\frac{\varepsilon}{24m^2},\\
  \label{eq:20}
  \dom(s_\xi)\cap \bigcup_{\alpha\in\varGamma_\xi}
  \bigcup_{n=p_\xi}^\infty B_\alpha^n=\emptyset.
\end{gather}
It clearly follows from~\eqref{eq:59} that there is a $t_\xi\supseteq s_\xi$ where
\begin{gather}
  \label{eq:3}
  \dom(t_\xi)=\bigcup_{\alpha\in\varGamma_\xi}\bigcup_{n<p_\xi}B_\alpha^n,\\
  \label{eq:13}
  \mu(x_\xi\cond [t_\xi])
  >\sqrt{1-\frac\varepsilon2}.
\end{gather}
By going to an uncountable subset $X\subseteq \oone$, 
we arrange that $p_\xi=p$ for all $\xi\in X$, and that $\mu([t_\xi])=\tau$ for all $\xi\in X$.   
Choose $l\ge p$ large enough so that 
\begin{equation}
  \label{eq:15}
  \sum_{n=l}^\infty h(n)<\frac{\tau^2\varepsilon}{24m^2}.
\end{equation}
And for each $\xi\in X$, put
\begin{equation}
  \label{eq:6}
  \vomega_\xi=\bigcup_{\alpha\in\varGamma_\xi}\bigcup_{n<l}B_\alpha^n.
\end{equation}

By going to an uncountable subsequence, we obtain $X'\subseteq X$ such that
$\{\vomega_\xi:\xi\in X'\}$ forms a $\Delta$-system, say with root
$\vomega$. Let $(\gamma^i:i<o)$ be the
strictly increasing enumeration of $\vomega$. By further
refinement, we can furthermore assume that 
\begin{multline}
  \label{eq:22}
  \bigl(\{(C^n_\alpha,D^n_\alpha:n<l):\alpha\in\varGamma_\xi\},t_\xi,(\gamma^i:i<o)\bigr)\\
  \cong\bigl(\{(C^n_\beta,D^n_\beta:n<l):\beta\in\varGamma_\eta\},t_\eta,(\gamma^i:i<o)\bigr)
\end{multline}
for all $\xi,\eta\in X'$, meaning that there is a finite partial injection $g$
on $\oone$ such that if every instance of each ordinal $\zeta$ appearing in the structure 
$\bigl(\{(C^n_\alpha,D^n_\alpha:n<l):\alpha\in\varGamma_\xi\},t_\xi,(\gamma^i:i<o)\bigr)$
is replaced with $g(\zeta)$ then 
$\bigl(\{(C^n_\beta,D^n_\beta: n<l):\linebreak 
\beta\in\varGamma_\eta\},t_\eta,(\gamma^i:i<o)\bigr)$ is obtained,
i.e.~$\{(C_\beta,D_\beta:n<l):\beta\in\varGamma_\eta\}=\bigl\{\bigl(\{s\circ g\inv:s\in
C_\alpha^n\},\{s\circ g\inv:s\in D_\alpha^n\}:n<l\bigr):\alpha\in\varGamma_\xi\bigr\}$,
$t_\eta=t_\xi\circ g\inv$ and $g(\gamma^i)=\gamma^i$ for all $i<o$.

\begin{claim}
\label{c-1}
$[t_\xi]\lands[t_\eta]\lands \bigl((y_\alpha^n\lands z_\beta^n)\lors
(y_\beta^n\lands z_\alpha^n)\bigr)=0$ 
for all $k\le n<p$, for all $\alpha\in\varGamma_\xi$ and all $\beta\in\varGamma_\eta$, 
for all $\xi,\eta\in X'$.
\end{claim}
\begin{proof}
Fix $\xi,\eta\in X'$ and $k\le n<p$.
For every $\alpha\in\varGamma_\xi$, since $\dom(t_\xi)\supseteq B_\alpha^n$, 
$[t_\xi]\lands y_\alpha^n$ is either 0 or $[t_\xi]$, and $[t_\xi]\lands z_\alpha^n$ is
either 0 or $[t_\xi]$. 
Similarly, for every $\beta\in\varGamma_\eta$, 
$[t_\eta]\lands y_\beta^n$ and $[t_\eta]\lands z_\beta^n$ are both either 0 or $[t_\eta]$. 
Now fix $\alpha\in\varGamma_\xi$ and $\beta\in\varGamma_\eta$.
The isomorphism witnessing~\eqref{eq:22} for $\xi$ and $\eta$ maps 
$(C^n_\delta,D^n_\delta)$ to $(C^n_\beta,D^n_\beta)$ for some $\delta\in\varGamma_\xi$, 
so that both
\begin{align} 
[t_\eta]\lands y_\beta^n=0&\Iff [t_\xi]\lands y_\delta^n=0,\\
[t_\eta]\lands z_\beta^n=0&\Iff [t_\xi]\lands z_\delta^n=0. 
\end{align}
However, by~\eqref{eq:39}, $x_\xi\lands y_\alpha^n\lands z_\delta^n=0$ and 
$x_\xi\lands y_\delta^n\lands z_\alpha^n=0$. And since $[t_\xi]\lands x_\xi\ne0$, this
implies that at most one of $[t_\xi]\lands y_\alpha^n$ and $[t_\xi]\lands z_\delta^n$
is nonzero, and at most one of $[t_\xi]\lands y_\delta^n$ and $[t_\xi]\lands
z_\alpha^n$ is nonzero.
It follows that both
$[t_\xi]\lands y_\alpha^n\lands[t_\eta]\lands z_\beta^n=0$ and 
$[t_\xi]\lands z_\alpha^n\lands[t_\eta]\lands y_\beta^n=0$, proving the claim.
\end{proof}

Note that $t_\xi$ is compatible with $t_\eta$ for all $\xi,\eta\in X'$
by~\eqref{eq:22}, since $t_\zeta\subseteq\vomega_\zeta$ for all $\zeta\in X$
by~\eqref{eq:3} and~\eqref{eq:6}.

\begin{claim}
\label{c-2}
$[t_\xi\cup t_\eta]$ and 
$\bigl\{\spr_{n=p}^{l-1}y^n_\alpha\lands z^n_\beta,
\spr_{n=p}^{l-1}y^n_\beta\lands z^n_\alpha\bigr\}$ 
are independently determined for all $\alpha\in\varGamma_\xi$ and
$\beta\in\varGamma_\eta$, for all $\xi\ne\eta$ in $X'$.
\end{claim}
\begin{proof}
Clearly $[t_\xi\cup t_\eta]$ is $\dom(t_\xi)\cup\dom(t_\eta)$-determined, and
$\bigl\{\spr_{n=p}^{l-1}y^n_\alpha\lands z^n_\beta,\linebreak
\spr_{n=p}^{l-1}y^n_\beta\lands z^n_\alpha\bigr\}$ 
is determined by coordinates in 
\begin{equation}
  \label{eq:58}
  \bigcup_{n=p}^{l-1}B^n_\alpha\cup\bigcup_{n=p}^{l-1}B^n_\beta.
\end{equation}
But $\dom(t_\xi)\cap\bigcup_{n=p}^{l-1}B^n_\beta\subseteq\dom(t_\xi)\cap\vomega_\eta$ 
is a subset of the root $\vomega$,
while the isomorphism from~\eqref{eq:22} fixes everything in the root, and thus
$\dom(t_\xi)\cap\bigcup_{n=p}^{l-1}B^n_\beta 
=\dom(t_\xi)\cap\bigcup_{n=p}^{l-1}B^n_\alpha=\emptyset$ by~\eqref{eq:3}. And similarly, 
$\dom(t_\eta)\cap\bigcup_{n=p}^{l-1}B^n_\alpha  
=\dom(t_\eta)\cap\bigcup_{n=p}^{l-1}B^n_\beta=\emptyset$. 
This establishes that $\dom(t_\xi)\cup\dom(t_\eta)$ is disjoint from the set
in~\eqref{eq:58}, as required.
\end{proof}

Now using~\eqref{eq:13} in Theorem~\ref{u-1}, 
we obtain $\xi\ne\eta$ in $X'$ such that
\begin{equation}
  \label{eq:38}
  \mu\bigl([t_\xi\cup t_\eta]\lands x_\xi\lands x_\eta\bigr)
  >\left(1-\frac\varepsilon2\right)\lands\mu([t_\xi\cup t_\eta]).
\end{equation}
And for all $\alpha\in\varGamma_\xi$ and $\beta\in\varGamma_\eta$,
\begin{equation}
  \label{eq:42}
  \begin{split}
  \mu\left([t_\xi\cup t_\eta]
    \lands\spr_{n=k}^\infty y_\alpha^n\lands z_\beta^n\right)
  &\le\mu\left([t_\xi\cup t_\eta]
    \lands\spr_{n=k}^{p-1} y_\alpha^n\lands z_\beta^n\right)\\
  &\hspace{-70pt}+\mu\left([t_\xi\cup t_\eta]
    \lands\spr_{n=p}^{l-1} y_\alpha^n\lands z_\beta^n\right)
  +\mu\left([t_\xi\cup t_\eta]
    \lands\spr_{n=l}^\infty y_\alpha^n\lands z_\beta^n\right)\\
  &\le0+\mu([t_\xi\cup t_\eta])
  \lands\mu\left(\spr_{n=p}^{l-1}y_\alpha^n\lands z_\beta^n\right)\\
  &\qquad
  +\mu\left(\spr_{n=l}^\infty y_\alpha^n\lands z_\beta^n\right)\\
  &\le \mu\bigl([t_\xi\cup t_\eta]\bigr)\cdot\sum_{n=p}^\infty
  3h(n)+\sum_{n=l}^\infty 3h(n)\\
  &<\frac{\varepsilon\cdot\mu([t_\xi\cup t_\eta])}{4m^2},
  \end{split}
\end{equation}
where Claims~\ref{c-1} and~\ref{c-2} are used for the second
inequality,~\eqref{eq:12} for the third,
and~\eqref{eq:21} and~\eqref{eq:15} are used for the fourth inequality; and
similarly
\begin{equation}
  \label{eq:45}
  \mu\left([t_\xi\cup t_\eta]
    \lands\spr_{n=k}^\infty y_\beta^n\lands z_\alpha^n\right)
  <\frac{\varepsilon\cdot\mu([t_\xi\cup t_\eta])}{4m^2}.
\end{equation}
Thus the measure of 
$[t_\xi\cup t_\eta]
\lands\spr_{\alpha\in\varGamma_\xi}\spr_{\beta\in\varGamma_\eta}
\spr_{n=k}^\infty (y_\alpha^n\lands z_\beta^n)
\lors(y_\beta^n\lands z_\alpha^n)$ is less
than~$\frac\varepsilon2\cdot\mu([t_\xi\cup t_\eta])$, which with~\eqref{eq:38} tells us that
\begin{multline}
  \label{eq:18}
  \mu\Biggl([t_\xi\cup t_\eta]\lands
    \Biggl(x_\xi\lands x_\eta\\
    -\spr_{\alpha\in\varGamma_\xi}\spr_{\beta\in\varGamma_\eta}
    \lbrak (\dot c_\alpha\cap\dot d_\beta\setminus k)
    \cup(\dot c_\beta\cap\dot d_\alpha\setminus k)
    \ne\emptyset\rbrak\Biggr)\Biggr)\\
    >(1-\varepsilon)\cdot\mu([t_\xi\cup t_\eta]),
\end{multline}
and in particular the condition~\eqref{eq:54} is satisfied.
This proves that $(\dot c_\alpha,\dot d_\alpha:\linebreak \alpha<\oone)$ is Souslin, thereby
completing the proof of Theorem~\ref{l-1}.

\section{$(\oone,\oone)$ Souslin gap}
\label{sec:destructible-gap}

We conclude the paper with a proof of Theorem~\ref{bigthm:2}, 
and thus Theorem~\ref{bigthm:1}, 
by constructing an $(\oone,\oone)$ Souslin gap in $L^0(\mu)\div\Fin$ for $\mu$ a
nonseparable measure. 
As we have seen, we can assume that $\mu=\mu_{(\theta)}$ for
some uncountable cardinal $\theta$.  Write $\random=\random_{(\theta)}$. 

Define a mapping $\phi:\theta\times\N\times\integers\to\Fin(\theta,\two)$  
where $\dom(\phi(\alpha,i,j))=[\alpha,\alpha+i)$ and the concatenation
\begin{equation}
\phi(\alpha,i,j)(\alpha+i-1)\phi(\alpha,i,j)(\alpha+i-2)\cdots \phi(\alpha,i,j)(\alpha)
\end{equation}
is the base $2$ representation of $j\mod{2^i}$. 
This can be expressed equivalently as
\begin{equation}
  \label{eq:4}
  \phi(\alpha,i,j)(\alpha+k)=\frac{j\mod{2^i}}{2^k}\land 1
  \espc\text{for all $k<i$.}
\end{equation}

Define $s_\alpha,t_\alpha\in \Fin(\theta,\two)^\N$ by
\begin{align}
  \label{eq:1}
  s_\alpha(n)
  &=\phi(\omega\cdot\alpha+n\lceil\log_2(n+1)\rceil,\lceil\log_2(n+2)\rceil,0),\\
  t_\alpha(n)
  &=\phi(\omega\cdot\alpha+n\lceil\log_2(n+1)\rceil,\lceil\log_2(n+2)\rceil,1).
\end{align}
Define random variables $\dot c_\alpha, \dot d_\alpha$ $(\alpha<\theta)$ by
\begin{align}
  \label{eq:5}
  \lbrak n\in \dot c_\alpha\rbrak&=[s_\alpha],\\
  \label{eq:2}
  \lbrak n\in \dot d_\alpha\rbrak&=[t_\alpha].
\end{align}
First of all, note that since $\log_2(n+2)\ge 1$ for all $n$, the third parameter in
$\phi$ ensures that
\begin{equation}
  \label{eq:28}
  \lbrak \dot c_\alpha\cap\dot d_\alpha
  =\emptyset\rbrak=1\espc\text{for all $\alpha$}.
\end{equation}
And for all $\alpha\ne\beta$, by stochastic independence,
$\lmeas n\in\dot c_\alpha\cap\dot d_\beta\rmeas
=2^{-2\lceil \log_2(n+2)\rceil}\linebreak\le\frac 1{(n+2)^2}$
for all $n$, and this sequence is summable which implies 
that\linebreak $[\dot c_\alpha]\lands[\dot d_\beta]=0$ in $L^0(\mu)\div\Fin$ 
(i.e.~$\dot c_\alpha\cap\dot d_\beta$ is finite with probability one).

On the other hand,
the length $\lceil\log_2(n+2)\rceil$ of the basic elements is short enough so that
$\lmeas n\in\dot c_\alpha\rmeas=\lmeas n\in\dot d_\alpha\rmeas
=2^{-\lceil\log_2(n+2)\rceil}$ and hence
\begin{equation}
  \label{eq:30}
  \frac1{2n+4}<\lmeas n\in\dot c_\alpha\rmeas\le\frac1{n+2},
\end{equation}
which is used to prove the following.

\begin{claim}
\label{c-3}
With probability one\textup: no subset of $\N$ interpolates the two families 
$\{\dot c_\alpha:\alpha<\theta\}$ and $\{\dot d_\alpha:\alpha<\theta\}$. In
particular, they form a gap in $L^0(\mu)\div\Fin$.
\end{claim}
\begin{proof}
Suppose that $r\in\cantorcube \theta$ is an $\random$-generic object (over $V$). 
Suppose that $e\in V[r]$ is a subset of $\N$ 
such that $\dot c_\alpha[r]\subseteqfnt e$ for all $\alpha$. 
Choose a countable $J\subseteq\theta$
large enough so that $e\in V[r\restriction I]$, 
where $I=\bigcup_{\alpha\in J}[\omega\cdot\alpha,\omega\cdot\alpha+\omega)$.
It follows from Kunen's Theorem~\ycite{Ku1}---stating that for any
$K\subseteq\theta$ in $V$, $r\in\cantorcube\theta$ is an
$\random_{(\theta)}$-generic object over $V$ iff $r\restriction K$ is
$\random_{(K)}$-generic over $V$ and $r\restriction\theta\setminus K$ is
$\random_{(\theta\setminus K)}$-generic over $V[r\restriction K]$---that 
$V[r]$ is a
forcing extension of $V[r\restriction I]$ by the measure
algebra $\Scal=\random_{({\theta\setminus I})}$ taken in $V[r\restriction I]$, 
and that for every $\alpha\notin J$, the $\Scal$-names for $\dot c_\alpha[r]$ and
$\dot d_\alpha[r]$, have the same definitions as given in~\eqref{eq:5}
and~\eqref{eq:2}, respectively.

In $V[r\restriction I]$: Fix $\alpha\notin J$. Let $\J$ be the set of all
$a\subseteq\N$ such that
\begin{equation}
  \label{eq:24}
  \lbrak a\cap\dot c_\alpha\text{ is finite}\rbrak\ne0.
\end{equation}
Since 
\begin{equation}
  \label{eq:14}
  (n+1)\lceil\log_2(n+2)\rceil-n\lceil\log_2(n+1)\rceil\ge\lceil\log_2(n+2)\rceil
\end{equation}
for all $n\in\N$, $\lbrak n\in\dot c_\alpha\rbrak$ ($n\in\N$) is a stochastically
independent sequence, and thus by Cauchy's criterion for infinite products,
\begin{equation}
\label{eq:35}
\begin{split}
a\in\J&\Iff \spr_{k=0}^\infty\ifm_{n\in a\setminus k}-\lbrak n\in\dot c_\alpha\rbrak\ne0\\
&\Iff\sup_{k\in\N}\prod_{n\in a\setminus k}1-\lmeas n\in\dot c_\alpha\rmeas=1\\
&\Iff\sum_{n\in a}\lmeas n\in\dot c_\alpha\rmeas<\infty\\
&\Iff\sum_{n\in a}(n+1)^{-1}<\infty\\
&\Iff a\in\ideal_{\frac1n}
\end{split} 
\end{equation}
where $\ideal_{\frac1n}$ is the well
known analytic ideal of all subsets of $\N$ on which the function $1\div(n+1)$ is
summable. In particular, $\J$ is a nonprincipal ideal, 
which means that $\lbrak\dot c_\alpha$ is infinite$\rbrak=1$. 

Since $\lbrak e\supseteqfnt\dot c_\alpha\rbrak\ne0$, $e\notin\J$.
However, $\J$ is also the ideal of all subsets $a\subseteq\N$ such that $\lbrak
a\cap\dot d_\alpha\text{ is finite}\rbrak\ne0$. Hence the proof is
complete, because we have shown that $e\cap \dot d_\alpha[r]$ is infinite.
\end{proof}

Notice that by~\eqref{eq:14}, for every $n\in\N$, 
$\bigcup_{\alpha<\oone}\bigl\{\lbrak n\in\dot c_\alpha\rbrak,
\lbrak n\in\dot d_\alpha\rbrak\bigr\}$ is determined by the coordinates
$\bigcup_{\alpha<\oone}[\omega\cdot\alpha+n\lceil \log_2(n+1)\rceil,
\omega\cdot\alpha+(n+1)\lceil\log_2(n+2)\rceil)$, and thusly the families are
independently determined for $m\ne n$, as condition~\eqref{eq:56} 
of Theorem~\ref{l-1} requires. Hence Theorem~\ref{l-1} entails that
$(\dot c_\alpha,\dot d_\alpha:\linebreak \alpha<\oone)$ is Souslin 
by letting $h\in\elspc1$ be given by $h(n)=\frac1{(n+2)^2}$, because 
\begin{equation}
  \label{eq:23}
  \lmeas n\in\dot c_\alpha\rmeas,\lmeas n\in\dot d_\alpha\rmeas
  =2^{-\lceil \log_2(n+2)\rceil}\le\frac1{n+2}\espc\text{for all $n$,}
\end{equation}
and therefore for all $\alpha\ne\beta$, by stochastic independence,
\begin{equation}
  \label{eq:26}
  \lmeas n\in\dot c_\alpha\cap\dot d_\beta\rmeas
  =\lmeas n\in\dot c_\alpha\rmeas\lmeas n\in\dot d_\beta\rmeas
  \le h(n),
\end{equation}
which with~\eqref{eq:28} assures condition~\eqref{eq:10}.

Now we build an $(\oone,\oone)$ Souslin gap $(\dot a_\alpha,\dot
b_\alpha:\alpha<\oone)$ where the $[\dot c_\alpha]\subseteqfntae[\dot a_\alpha]$ and $[\dot
d_\alpha]\subseteqfntae[\dot b_\alpha]$ in $L^0(\mu)\div\Fin$ for all
$\alpha<\oone$. Therefore, since $(\dot c_\alpha,\dot d_\alpha:\alpha<\oone)$ forms
a gap by Claim~\ref{c-3}, $(\dot a_\alpha,\dot b_\alpha:\alpha<\oone)$ will automatically
be a gap. The idea is to limit the augmentation of the measures of $\lbrak n\in\dot
c_\alpha\rbrak$ and $\lbrak n\in\dot d_\alpha\rbrak$ so that Theorem~\ref{l-1} still
applies. 

The $(\oone,\oone)$ gap we have been striving towards is obtained 
from an `ascending tower' of sorts (espc.~conditions~\eqref{item:7} and~\eqref{item:8} 
below). Namely, a sequence $T_\alpha\in\Fin(\oone)^\N$ ($\alpha<\oone$), 
where $\Fin(\oone)$ denotes the collection of all finite subsets of $\oone$,
satisfying:
\begin{enumerate}[(m)]
\item\label{item:7} $\alpha\in T_\alpha(n)$ for all $n$,
\item\label{item:8} $T_\xi(n)\subseteq T_\alpha(n)$ for all but finitely many $n$, for all
  $\xi<\alpha$,
\item\label{item:6} $|T_\alpha(n)|\le\sqrt[3]{n+1}$ for all $n$.
\save
\end{enumerate}
Such a tower can easily be constructed by recursion on $\alpha$, 
by adding the requirement
\begin{enumerate}[(m)]
\restore
\item $\displaystyle\lim_{n\to\infty}\frac{|T_\alpha(n)|}{\sqrt[3]{n+1}}=0$
\end{enumerate}
to carry the recursion through.

Now $(\dot a_\alpha,\dot b_\alpha:\alpha<\oone)$ is defined by
\begin{equation}
  \label{eq:41}
  \lbrak n\in\dot a_\alpha\rbrak=\spr_{\xi\in T_\alpha(n)}[s_\xi(n)]
  \AND\lbrak n\in\dot b_\alpha\rbrak=\spr_{\xi\in T_\alpha(n)}[t_\xi(n)]
\end{equation}
for all $n\in\N$, for all $\alpha<\oone$. Then in fact, $\lbrak\dot
a_\alpha\supseteq\dot c_\alpha\rbrak=1$ and $\lbrak\dot b_\alpha\supseteq\dot
d_\alpha\rbrak=1$ by~\eqref{item:7}, and with~\eqref{item:8} it follows that both
$\{\dot a_\alpha:\alpha<\oone\}$ and $\{\dot b_\alpha:\alpha<\oone\}$ have order
type $\oone$. Applying Lemma~\ref{l-5} to the
independently determined family $\bigl\{\{[s_\xi(n)],[t_\xi(n)]\}:\xi<\oone\bigr\}$, with
both of the subsets $A$ and $B$ equal to $T_\alpha(n)$, and using~\eqref{eq:23} with
condition~\eqref{item:6} yields
\begin{equation}
  \label{eq:62}
  \begin{split}
  \lmeas n\in\dot a_\alpha\cap\dot b_\beta\rmeas&\le\lmeas n\in\dot a_\alpha\rmeas
  \lmeas n\in\dot b_\beta\rmeas\\
  &\le\left(\frac{\sqrt[3]{n+1}}{n+2}\right)^2\\
  &<(n+1)^{-\frac 43}
  \end{split}
\end{equation}
for all $n\in\N$, for all $\alpha,\beta<\oone$. Since this extension of the original
nonlinear gap still satisfies~\eqref{eq:56}, Theorem~\ref{l-1} applies with the function in
$\elspc1$ given by\footnote{In this construction, any function which goes to infinity
  and is everywhere $\ge 1$ can be
  used in place of $\sqrt[3]{n+1}$ as long as the resulting function corresponding
  to~\eqref{eq:66} is in $\elspc1$.}
\begin{equation}
  \label{eq:66}
  n\mapsto (n+1)^{-\frac 43}.
\end{equation}
Therefore $(\dot a_\alpha,\dot b_\alpha:\alpha<\oone)$ 
is an $(\oone,\oone)$
Souslin gap in the Boolean algebra $(L^0(\mu)\div\Fin,\subseteqfntae)$, 
concluding the paper.



\medskip
\indent\textsc{\small Graduate School of Science and Technology, 
Kobe University, Japan}

\indent\small{\textit{E-mail address}: \texttt{j\_hirschorn@yahoo.com}}

\indent\small{\textit{URL}: \texttt{http://www.logic.univie.ac.at/\textasciitilde hirschor/}}

\end{document}